\documentclass[12pt,reqno]{amsart}
\usepackage{graphicx}
\usepackage{hyperref}
\usepackage{psfrag}
\usepackage{pxfonts}
\usepackage{mathrsfs}
\usepackage{color}
\usepackage{xcolor}
\usepackage{amsmath,amssymb}
\newcommand\numberthis{\addtocounter{equation}{1}\tag{\theequation}}
\numberwithin{equation}{section}
\newcommand{\diff}{\operatorname{Diff}}

\theoremstyle{plain}
\newtheorem{maintheorem}{Theorem}

\newcommand{\R}{\mathbb{R}}
\newcommand{\N}{\mathbb{N}}
\newcommand{\Z}{\mathbb{Z}}

\newcommand{\T}{\mathbb{T}}
\newcommand{\F}{\mathcal{F}}
\newcommand{\vol}{\mathrm{Vol}}
\newcommand{\jac}{\mathrm{Jac}}
\newcommand{\dime}{\mathrm{dim}}

\newtheorem{theorem}{Theorem}[section]

\newtheorem{corollary}[theorem]{Corollary}
\newtheorem{proposition}[theorem]{Proposition}

\newtheorem{definition}[theorem]{Definition}
\newtheorem{Conjecture}[theorem]{Conjecture}

\theoremstyle{remark}

\newtheorem{problem}{Problem}

\begin{document}

\thanks{}

\author{J. Santana Costa}
\address{DEMAT-UFMA S\~{a}o Luis-MA, Brazil.}
\email{jsc.costa@ufma.br}
\thanks{J. S. Costa was supported by CAPES-PROEX and CNPq process 141224/2013-4}

%\author{F. Micena}
%\address{
%  IMC-UNIFEI Itajub\'{a}-MG, Brazil.}
%\email{fpmicena@gmail.com}

\author{A. Tahzibi}
\address{Departamento de Matem\'atica,
  ICMC-USP S\~{a}o Carlos-SP, Brazil.}
\email{tahzibi@icmc.usp.br}
\thanks{A. Tahzibi was supported by the FAPESP thematic project 2017/06463-3 and the CNPq productivity fellowship.}

\renewcommand{\subjclassname}{\textup{2000} Mathematics Subject Classification}

\date{\today}

\setcounter{tocdepth}{2}

\title[Rigidity of Lyapunov exponents ]{Rigidity of Lyapunov exponents for derived from Anosov diffeomorphisms}
\maketitle
\begin{abstract}
 For a class of volume preserving partially hyperbolic diffeomorphisms (or non-uniformly Anosov) $f\colon {\T}^d\rightarrow{\T}^d$ homotopic to linear Anosov automorphism, we show that the sum of the positive (negative) Lyapunov exponents of $f$  is bounded above (resp. below) by the sum of the positive (resp. negative) Lyapunov exponents of its linearization. We show this for some classes of derived from Anosov (DA) and non-uniformly hyperbolic systems with dominated splitting, in particular for examples described by C. Bonatti and \mbox{ M. Viana \cite{bonattiviana2000srb}}.  The results in this paper address a flexibility program by J. Bochi, A. Katok and F. Rodriguez Hertz \cite{BKR}.

%Our results generalize results of (Micena, Tahzibi,13), (Micena, Tahzibi,14) and (Saghin, Xia).
\end{abstract}

%----------------------------------------------------------
\section{Introduction}
%----------------------------------------------------------

The Lyapunov exponents come from the study of differential equations in the thesis of A. M. Lyapunov \cite{lyapunov1992general}. They were systematically introduced into ergodic theory by works of H. Furstenberg and H. Kesten \cite{furstenberg1960products} and V. Oseledets \cite{oseledec1968multiplicative}.
%When Oseledets was PhD student, he was studying metric entropy of dynamical systems, a concept developed by A. N. Kolmogorov and Ya. G. Sinai (Oseledets' Advisor) and it became clear to him that positive entropy was directly related to the expansion rates of the system. Positive responses to these ideas were given by
They are directly related to the expansion rates of the system and also to the positive metric entropy, using for example
%by D. Ruelle and  Ya.  Pesin, known as 
the Margulis-Ruelle inequality and the Pesin entropy formula.  Ya. Pesin explored a lot the concept of Lyapunov exponents developing a rich theory of non-uniformly hyperbolic systems, which are the systems with non-zero Lyapunov exponents.

However in general these systems are not robust: they do not form an open set. In the works \cite{bochi2002genericity} and \cite{bochi2005lyapunov}, the authors show that if a  non-uniformly hyperbolic system does not admit dominated splitting, it can be approximated by systems with zero Lyapunov exponents in $C^1$ topology. Usually, the presence of dominated decomposition guarantees higher regularity in  $C^1$ topology, see for example \cite{saghin2020regularity}. In the non-invertible setting,  
\cite{andersson2022non} showed the existence of $C^1$ open sets without dominated decomposition such that integrated Lyapunov exponents vary continuously with the dynamics in the $C^1$ topology. 

 For a linear Anosov automorphism of the torus $A \colon  {\T}^d\rightarrow{\T}^d$ the Lyapunov exponents are constant 
 and indeed they are equal to the logarithm of the norm of eigenvalues of $A$. In general, Lyapunov exponents 
cannot be explicitly calculated. The regularity of Lyapunov exponents with respect to dynamics and invariant measures  is a subtle question. Let us recall a flexibility conjecture by J. Bochi, A. Katok and F. Rodriguez Hertz. Let $\diff_m^{\infty} (M)$ be the set of $m-$preserving diffeomorphisms (volume preserving)
\mbox{$f : M \rightarrow M$} of class $C^{\infty}$.

\begin{Conjecture}[\cite{BKR}]    Given a connected component $C \subset \diff_m^{\infty} (M)$ and any list of numbers $\xi_1 \geq \cdots \geq \xi_d$ with $\sum_{i=1}^{d} \xi_i =0,$   there exists an ergodic diffeomorphism $f \in C$ such that $\xi_i, i=1, \cdots, d$ are the Lyapunov exponents with respect to $m.$
\end{Conjecture}

Moreover, in the setting of volume preserving Anosov diffeomorphisms, they posed the following problem:

\begin{problem}[\cite{BKR}] \label{strong} (Strong flexibility)
Let $A \in SL(d, \mathbb{Z})$ be a hyperbolic linear transformation inducing conservative Anosov diffeomorphism $F_A$ on $\mathbb{T}^d$ with Lyapunov exponents: $\lambda_1 \geq \lambda_2 \cdots \geq \lambda_{u} > 0 > \lambda_{u+1} \geq \cdots \geq \lambda_d.$ Given any list of numbers $\xi_1 \geq \xi_2 \geq \cdots \geq \xi_u \geq \xi_{u+1} \geq \cdots \geq \xi_d$ such that 
\begin{enumerate}
\item $\sum_{i=1}^{d} \xi_i = 0,$
\item $\sum_{i=1}^{u} \xi_i \leq \sum_{i=1}^{u} \lambda_i,$
\end{enumerate}
does there exist a conservative Anosov diffeomorphism $f$ homotopic to $F_A$ such that $\{\xi_i\}$ is the list of all Lyapunov exponents with respect to volume measure?

\end{problem}

In this work, we study a subset of transformations $f\colon {\T}^d\rightarrow{\T}^d$ homotopic to a  linear Anosov automorphism $A\colon {\T}^d\rightarrow{\T}^d$ where $f$ has some hyperbolicity (partial hyperbolicity or  non-uniform hyperbolicity). See section 2 for the definitions. In particular, we show that for a class of partially hyperbolic diffeomorphisms and homotopic to Anosov linear automorphism, not all lists of numbers can be realized as Lyapunov exponents. More precisely, in such a class of dynamics, we need the conditions (1) and (2) imposed in the Problem \ref{strong}.

 This type of result also appear in   \cite{BKR}, \cite{CostaMicena2017pathological}, \cite{micena2017new}, \cite{micena2013regularity}, \cite{micena2016unstable}
and \cite{SaghinXia2009geometric}. In \cite{micena2013regularity} it has been proved that for any conservative partially hyperbolic systems in the torus ${\T}^3$ the stable or unstable Lyapunov exponent is bounded by the stable or unstable Lyapunov exponent of its linearization.  It is well worth mentioning that Carrasco-Saghin \cite{CS} constructed a $C^\infty$ and volume preserving example on $\mathbb{T}^3$ that shows that the largest Lyapunov exponent of a diffeomorphism in the homotopy class of an Anosov linear automorphism of $\mathbb{T}^3$ may be larger than the largest Lyapunov exponent of $A$. However, in their example, $f$ does not admit a three-bundle partially hyperbolic splitting. 

 We also mention that our results address (and give a negative answer in some special cases of derived from Anosov diffeomorhisms) a question in \cite{CS}:
Does there exist a derived from Anosov diffeomorphism $f$ such that the sum of the positive Lyapunov exponents of $f$ is larger than the sum of the positive Lyapunov exponents of its linear part?

%From \cite{SaghinXia2009geometric} it can be verified that the same holds for conservative partially hyperbolic in the torus ${\T}^d$ ones that are close to linear Anosov.

%\newpage

\section{Definitions and Statements of Results}

\begin{definition}
Let $M$ be a closed manifold.
A diffeomorphism $f\colon  M \rightarrow M$ is called a partially hyperbolic diffeomorphism if there is a suitable norm $||\cdot||$ and the tangent bundle $TM$ admits a $Df$-invariant decomposition $TM =  E^s \oplus E^c \oplus E^u$ such that for all unitary vectors $v^{\sigma} \in E^{\sigma }_x, \sigma \in \{s,c,u\}$ and every $x \in M$ we have:

$$ ||D_x f v^s || < ||D_x f v^c || < ||D_x f v^u ||,$$
moreover,

$$||D_x f v^s || < 1 \;\mbox{and}\; ||D_x f v^u || > 1 .$$

\end{definition}

Every diffeomorphism of the torus $\mathbb{T}^d$ induces an automorphism of the fundamental group and there exists a unique linear
diffeomorphism $f_{\ast}$ which induces the same automorphism on $\pi_1(\mathbb{T}^d).$ The diffeomorphism $f_{\ast}$ is called linearization of $f.$

 Anosov diffeomorphisms can be considered partially hyperbolic systems with $E^c = 0.$  However, in this paper whenever we consider an Anosov diffeomorphism as a partially hyperbolic system, we mean that there exists a non-trivial (center bundle) partially hyperbolic decomposition. 
 
 Clearly, a partially hyperbolic diffeomorphism may have various partially hyperbolic invariant decompositions. We will consider several (non-necessarily disjoint) categories of partially hyperbolic diffeomorphisms indexed by the dimensions of invariant bundles. 
 
 More precisely, we say $f \in Ph_{(d_s, d_u)}(M)$ if $f$ admits a partially hyperbolic decomposition with $\dim(E^{\sigma}) = d_{\sigma},$ for  $ \sigma \in \{s, u\}.$ Clearly $\dim(E^c) = \dim(M) - d_s - d_u.$
For instance, consider the cat map $A$ on $\T^2$ induced by matrix $
\left(
  \begin{array}{cc}
    2 & 1 \\
    1 & 1 \\
  \end{array}
\right)$, then $A \times A$ is a partially hyperbolic (in fact Anosov) diffeomorphism and belongs to $Ph_{(1, 1)} \cap Ph_{(2,1)} \cap Ph_{(1, 2)}.$

\begin{definition} Let $f\colon  \mathbb{T}^d \rightarrow \mathbb{T}^d $ be a partially hyperbolic diffeomorphism, $f$ is called a
derived from Anosov (DA) diffeomorphism if its linearization $f_{\ast}\colon  \mathbb{T}^d \rightarrow \mathbb{T}^d $ is a linear Anosov automorphism.
\end{definition}

\begin{definition}
We say that the diffeomorphism $f\colon  M\rightarrow M$ admits a dominated splitting if there is an invariant (by $Df$) continuous decomposition 
$TM=E\oplus F$ and constants  $0<\nu<1, C >0$ such that
$$
\frac{||Df^n|_{E(x)}||}{||Df^{-n}|_{F(f^n(n))}||^{-1}}\leq C\nu^n,\,\,\forall x\in M, n>0.
$$
\end{definition}

 Let us recall a simple  version of the Oseledets' theorem.

\begin{theorem} \label{oseledets} \cite{oseledec1968multiplicative}
Let $f\colon  M\rightarrow M$ be a $C^1$ diffeomorphism, then there is a full probability borelian set $\mathcal{R}$
(this is, $\mu(\mathcal{R})=1$ for all $f-$invariant probability measure $\mu$) such that for each $x\in \mathcal{R}$ there is a decomposition $T_xM=E_1\oplus \cdots E_{k(x)}$ and constants $\lambda_1,\ldots,\lambda_{k(x)}$ such that
$$
\displaystyle\lim_{n\rightarrow \infty}\frac{1}{n}||Df_{x}^{n}(v)||=\lambda_i,
$$ 
for all $v\in E_i$. The $\lambda_i(x)$ is called a Lyapunov exponent. Moreover, for any $ 1 \leq j \leq k(x) $
$$  \displaystyle\lim_{n\rightarrow \infty}\frac{1}{n} det Df_{x}^{n}|_{F_j(x)}= \sum_{i=1}^{j} d_i \lambda_i,$$ where $F_j = E_1 \oplus \cdots \oplus E_j$ and $d_i= dim (E_i).$
\end{theorem}
In this paper, when $m$ refers to a measure, it stands for a fixed  Lebesgue measure on $\mathbb{T}^d$.
\subsection{Lyapunov exponents of partially hyperbolic diffeomorphisms}

In our first main theorem, we prove that linear Anosov diffeomorphisms maximize (resp. minimize) the sum of unstable  (resp. stable) Lyapunov exponents, in any homotopy path totally inside partially hyperbolic diffeomorphisms. 

The notion of metric entropy appears naturally when dealing with Lyapunov exponents. For a partially hyperbolic diffeomorphism, it is possible to define entropy of an invariant measure along unstable foliation and there is variational principle results for such entropy. A $u-$maximal entropy measure is a measure which attains the supremum of unstable entropy among all invariant measures. See Subsection \ref{maximizing} and references for more details.
\begin{maintheorem}\label{teor B}
Let $f\colon {\T}^d\rightarrow {\T}^d$ be a $C^2$ volume preserving DA diffeomorphism in $Ph_{d_s, d_u}$  with  linearization $A\colon {\T}^d\rightarrow {\T}^d$  such that
 $f$ and $A$ are homotopic by a path fully contained in $Ph_{d_s, d_u}({\T}^d)$,
then
$$
\sum_{i = 1}^{d_u} \lambda^u_i(f, x) \leq  \sum_{i = 1}^{d_u} \lambda^u_i(A)
\,\,\,\,
\mbox{and}
\,\,\,\,
\sum_{i = 1}^{d_s} \lambda^s_i(f, x) \geq  \sum_{i = 1}^{d_s} \lambda^s_i(A),
$$
for $m-a.e.$  $x \in {\T}^d$ and,
the first (resp. second) inequality is strict unless $m$ is a measure of $u-$maximal entropy (resp. $u-$maximal for $f^{-1}$). 
\end{maintheorem}

This is related to the  result of \cite{micena2013regularity}: Let $f\colon  \mathbb{T}^3 \rightarrow \mathbb{T}^3$ be a $C^2$ conservative derived from Anosov partially hyperbolic diffeomorphism with linearization $A$. Then, $\lambda^u(x) \leq \lambda^u(A)$ for almost every $x.$  We recall that the authors used quasi-isometric property of unstable foliation for the three dimensional derived from Anosov diffeomorphisms. They do not assume that $f$ is in the same connected component of $A$.

\begin{maintheorem}\label{teor A}
Let $f\colon {\T}^d\rightarrow {\T}^d$ be a $C^2$ volume preserving DA diffeomorphism in $Ph_{d_s, d_u}$  with  linearization $A\colon {\T}^d\rightarrow {\T}^d$ belonging to $Ph_{d_s, d_u}.$
Suppose that for $\sigma \in \{s, u\}$
 there is $(d - d_{\sigma})-$dimensional subspace $P_{\sigma} \subset \mathbb{R}^d, $ such that $\angle (E^{\sigma}_f(x), P_{\sigma}) > \alpha > 0,$ for all $x \in \mathbb{R}^d,$ ($E^{\sigma}_f$ stands for the lift of the bundles)
then
$$
\sum_{i = 1}^{d_u} \lambda^u_i(f, x) \leq  \sum_{i = 1}^{d_u} \lambda^u_i(A)
\,\,\,\,
\mbox{and}
\,\,\,\,
\sum_{i = 1}^{d_s} \lambda^s_i(f, x) \geq  \sum_{i = 1}^{d_s} \lambda^s_i(A),
$$
for $m-a.e.$  $x \in \mathbb{T}^d.$
\end{maintheorem}

 Observe that the assumption on the existence of $P_{\sigma}$ in the above theorem is a mild condition and is satisfied if the stable and unstable bundles of $f$ do not vary too much. Indeed, we require that $\{ E^{\sigma}_f(x), x \in \mathbb{R}^d\}$ is not dense in the Grassmannian of $d_{\sigma}$ dimensional subspaces.

\subsection{Lyapunov exponents of  non-uniformly Anosov systems}\label{nonuniformanosov}
In this section, we deal with dynamics that are not necessarily partially hyperbolic. 
However,  they enjoy some dominated splitting property: Oseledets' splitting (stable/unstable) is dominated.

\begin{definition}
A $C^2$- volume preserving diffeomorphism $f$ is called non-uniformly Anosov if it admits an $f-$invariant dominated decomposition $TM= E \oplus F$  such that $\lambda^{F}_i(x)>0$ (all Lyapunov exponents along $F$) and $\lambda^{E}_i(x)<0$ (all Lyapunov exponents along $E$) for $\mu-$a.e $x$.
\end{definition}

%\begin{question}
%Is it true that any $C^2$-volume preserving non-uniformly Anosov diffeomorphism is ergodic if it is topologically transitive?
%\end{question}

In the next theorem, we compare Lyapunov exponents of a non-uniformly Anosov diffeomorphism with those of linear Anosov automorphism if they are homotopic. Observe that by definition, the number of negative (positive) Lyapunov exponents of a non-uniformly Anosov diffeomorphism is constant almost everywhere and coincides with the dimension of the bundles in the dominated splitting.  However, the Lyapunov exponents may depend to the orbit, as we do not assume ergodicity.
In fact, it is interesting to know whether in general all topologically transitive non-uniformly Anosov diffeomorphism are ergodic or not.  We conjecture that all topologically transitive non-uniformly Anosov diffeomorphisms are ergodic.

Let $A\colon  {\T}^d\rightarrow {\T}^d$ be a linear Anosov diffeomorphism such that $T \mathbb{T}^d = E^s_{A} \oplus E^u_{A}.$ Denote by \mbox{$d_s =\dime  (E^s_A)$} and  $d_u =\dime(E^u_A).$ After changing to an equivalent Riemannian norm, if necessary, there are $0<\lambda<1<\gamma$ such that $\|A|_{E^s}\| \leq \lambda$, $m(A|_{E^u}) \geq \gamma$ and $E^s_A$ is orthogonal to $E^u_A.$  Indeed, to use fewer constants in the proofs, we assume that $E^u_i$ and $E^u_j, i \neq j$ are orthogonal, where $E^u_i$'s are generalized eigenspaces of $A$ which coincide with the Oseledets decomposition.  Recall that for a linear transformation $T$,  $m(T):= \|T^{-1}\|^{-1}.$ We refer to $\gamma, \lambda$ as rates of hyperbolicity of $A$. Observe that Lyapunov exponents of any diffeomorphism are independent of the choice of the equivalent norm.

\begin{maintheorem}\label{TeoDecmDominada}
	Let $f\colon  {\T}^d\rightarrow {\T}^d$ be a $C^{2}$ conservative non-uniformly Anosov diffeomorphism with 
  dominated decomposition $TM=E\oplus F$ , homotopic to $A$ such that
\begin{enumerate}
  \item $\dime (E) =d_s$\,\, and \,\, $\dime (F) =d_u,$
  \item $E(x)\cap E^u_A=\{0\}$ and $F(x)\cap E^s_A=\{0\}$
  \item $||Df|_E||<\gamma$ and $m(Df|_F) >\lambda$,
\end{enumerate}
		 Suppose  that the distributions $E$ and $F$ are integrable, then $\displaystyle\sum_{i=1}^{d_u}\lambda^{F}_{i}(f,x)\leq\displaystyle \sum_{i=1}^{d_u}\lambda_i^u(A)$ and
	$\displaystyle\sum_{i=1}^{d_s}\lambda^{E}_i(f,x)\geq\displaystyle\sum_{i=1}^{d_s}
	\lambda_i^s(A)$ for Lebesgue almost every $x\in {\T}^d$.
\end{maintheorem}

Let us comment on the hypotheses: The first and second one ask for some compatibility of invariant bundles and the third one asks that any possible expansion in the dominated bundle $E$ is less  than the expansion rate of $A$ and similarly any possible contraction in the dominating bundle $F$ is weaker than the contraction rate of $A.$

C. Bonatti and M. Viana \cite{bonattiviana2000srb} constructed the first examples of robustly transitive diffeomorphisms which are not partially hyperbolic, which was later generalized in higher dimensions by \cite{tahzibi2004stably}. Those classes of examples satisfy the hypotheses of the above theorem.

\begin{theorem}[\cite{bonattiviana2000srb}, \cite{buzzifischer2013entropic}, \cite{tahzibi2004stably}]\label{Teo Bv}
	There is an open set $\mathcal{U}\subset \diff^1_m({\T}^d)$ such that any $C^2$-diffeomorphism in  $\mathcal{U}$ satisfies all hypotheses of Theorem C and the bundles $E$ and $F$ are integrable.
\end{theorem}

In the above theorem, the integrability of invariant subbundles $E$ and $F$ is proved in \cite{buzzifischer2013entropic}.

%\marginpar{vou escrever com mais cuidado no enunciado deste dos chilenos}
%Related to this Theorem, in \cite{saghin2020regularity} have a strong result about the regularity of the Lyapunov exponents. They show that for families of diffeomorphisms with dominated splitting, the map that takes $f$ in the sum of its Lyapunov exponents is differentiable.  

%By theorems \ref{TeoDecmDominada} and \ref{Teo Bv} we obtain the following corollary.

%\begin{corollary}
%There is an open set $\mathcal{U}\subset \diff^1_m({\T}^4)$ where each $f\in \mathcal{U}$ is a $C^{2}$ conservative non uniformly Anosov diffeomorphism with dominated decomposition $TM=E\oplus F$, homotopic to a linear Anosov diffeomorphism $A$ such that
%$\displaystyle\sum_i\lambda_i^{F}(f,x)\leq\displaystyle \sum_i\lambda_i^u(A)$ and
%	$\displaystyle\sum_i\lambda_i^{E}(f,x)\geq\displaystyle\sum_i
%	\lambda_i^s(A)$ for Lebesgue almost everywhere $x\in {\T}^4$.
%\end{corollary}

\subsection{Regularity of Foliations}

 Any partially hyperbolic diffeomorphism admits invariant stable and unstable foliations. As the proofs of our main results are based on the regularity of such foliations we need to recall some basic definitions and results in this subsection.
Consider ${\F}$ a foliation in $M$, $B$ a  ${\F}-$foliated box   and $m$ the Lebesgue measure in $M$. Denote by $\vol_{L_x}$ the Lebesgue measure of leaf
$L_x$ and $m_{L_x}$  the disintegration of the 
Lebesgue measure $m$ along the leaf $L_x$. The disintegrated measures are in fact projective class of measures. However whenever we fix a compact foliated box $B$ then we can use the Rohlin disintegration theorem for the normalized restriction of $m$ on $B$ to get probability conditional measures. So in what follows, after fixing a foliated box $B$ by $m_{L_x}$ and $\vol_{L_x}$ we understand probability measures whose support is inside the plaque of $L_x$ which contains $x \in B$ and is inside $B$.

\begin{definition}
We say that the foliation ${\F}$ is \textit{upper leaf wise absolutely continuous} if for any foliated box $B$,  $m_{L_x}\ll\vol_{L_x}$
for $m$-almost everywhere $x\in B$. Equivalently, if given a set $Z\subset B$ such that
$\vol_{L_x}(Z\cap L_x)=0$
for $m$-almost everywhere $x \in B$, then $m(Z)=0$.
\end{definition}

\begin{definition}
We say that the foliation ${\F}$ is \textit{lower leafwise absolutely continuous} 
if $\vol_{L_x}\ll m_{L_x}$
for $m$-almost everywhere $x\in B$. Equivalently, if given a set $Z\subset B$ such that $m(Z)=0$,
then $\vol_{L_x}(Z\cap L_x)=0$
for $m$-almost everywhere $x \in B$.
\end{definition}

\begin{definition}
We say that the foliation ${\F}$ is \textit{leafwise absolutely continuous} if \linebreak $\vol_{L_x}\sim m_{L_x}$ 
(this is, $m_{L_x}\ll\vol_{L_x}$ and
 $\vol_{L_x}\ll m_{L_x}$)
for $m$-almost everywhere $x\in B$.
\end{definition}

\begin{proposition}\cite{brinpesin}
If $f\colon  M\rightarrow M$ is a $C^2$ partially hyperbolic diffeomorphism, then the foliations $W^s$ and $W^u$ are leafwise absolutely continuous.
\end{proposition}

 We also mention the result of Y. Pesin for non-uniformly hyperbolic systems, see (Theorem 4.3.1, \cite{BarreiraPesin2002lyapunov}), which shows absolute continuity of 
 local Pesin's laminations.
 
\section{Proof of Results}

The theorems  \ref{teor B} and \ref{teor A} are obtained by the following result.

\begin{theorem}\label{Teo formas}
Let $f\colon {\T}^d\rightarrow {\T}^d$ be a $C^2$ volume preserving partially hyperbolic diffeomorphism such that there are  closed non-degenerate $d_u-$forms and $d_s-$ forms on $E^u$ and $E^s$, respectively. Suppose that $A$, the linearization of $f$,
is partially hyperbolic and $\dime E^{\sigma}_f=\dime E^{\sigma}_A$, $\sigma\in \{s,c,u\}$ then
$$
\displaystyle\sum_{i=1}^{d_u}\lambda_i^u(f,x)\leq\displaystyle \sum_{i=1}^{d_u}\lambda_i^u(A)
\,\,\,\,\mbox{and}\,\,\,\,
\displaystyle\sum_{i=1}^{d_s}\lambda_i^s(f,x)\geq\displaystyle\sum_{i=1}^{d_s}\lambda_i^s(A)
$$
for $m-a.e.$   $x\in {\T}^d$.% onde $\mathcal{R}$ é o conjunto regular de Lyapunov.
 \end{theorem}

For the proof of this theorem we use a result of R. Saghin, \cite{Saghin2014}.  Let $f\colon M\rightarrow M$ be a diffeomorphism and  $W$ a $f-$invariant foliation on $M$, that is, $f(W(x))=W(f(x))$  and $B_r(x,f)$ be the ball of the leaf
 $W(x)$ with radius
 $r$ centered at $x$. We say that
$$
\chi_W(f,x)=\displaystyle\limsup_{n\rightarrow\infty}\frac{1}{n}\log \vol(f^n(B_r(x,f)))
$$
is the \emph{volume growth rate of the foliation at} $x$ and
$$
\chi_W(f)=\displaystyle\sup_{x\in M}\chi_W(f,x)
$$
is the \emph{volume growth rate} of $W$. Here Vol stands for the $\vol_{W}$ which is the induced volume to the leaves of $W$. We use this notation along the paper except when it may create confusion.

When $f$ is a partially hyperbolic diffeomorphism we denote by $\chi_u(f)$ the volume growth of unstable foliation $W^u$.

\begin{theorem}[\cite{Saghin2014}]\label{Teo Saghin}
	Let $f\colon  M\rightarrow M$ be a $C^1$  partially hyperbolic diffeomorphism such that there is a closed
	 non-degenerate $d_u-$form on the unstable bundle $E^u$, then $\chi_u(f)=\log sp(f_{*,d_u})$, where $f_{*,d_u}$ is the induced map in the
	$d_u$-cohomology of De Rham.
\end{theorem}

\begin{proposition}\label{Prop vol formas}
	Let $f\colon {\T}^d\rightarrow {\T}^d$ be a $C^1$  partially hyperbolic diffeomorphism 
 admitting a closed non-degenerate $d_u-$form on the unstable bundle $E^u$. For fixed 
 $x\in {\T}^d, r >0$ the balls $B_r(x,f)\subset W^u(x,f)$ and $B_r(x,A)\subset W^u(x,A)$ 
 satisfy the following$\colon$ 
 Given $\varepsilon>0$, there is $n_0 \in \mathbb{N}$ such that if $n>n_0$ we have
	$$
	 \vol_{W^u(f)}(f^n(B_r(x,f)))\leq(1+\varepsilon)^n\vol_{W^u(A)}(A^n(B_r(x,A))).
	$$
\end{proposition}

\begin{proof}
	By Theorem \ref{Teo Saghin} we have that
	$$
	\chi_u(f,x)\leq\chi_u(f)=\log sp(f_{*,u})=\log sp(A_{*,u})=\chi_u(A)=\chi_u(A,x),
	$$
	that is
	$$
	\displaystyle\limsup_{n\rightarrow\infty}\frac{1}{n}\log \vol_{W^u(f)}(f^n(B_r(x,f)))\leq \displaystyle\lim_{n\rightarrow\infty}\frac{1}{n}\log  \vol_{W^u(A)}(A^n(B_r(x,A))),
	$$
	therefore give $\varepsilon>0$, there is $n_0$ such that if $n>n_0$,
	$$
	\frac{1}{n}\log \vol_{W^u(f)}(f^n(B_r(x,f)))\leq \frac{1}{n}\log \vol_{W^u(A)}(A^n(B_r(x,A)))+\frac{1}{n}\log (1+\varepsilon)^n.
	$$
	Then
	$$
	\vol_{W^u(f)}(f^n(B_r(x,f)))\leq(1+\varepsilon)^n\vol_{W^u(A)}(A^n(B_r(x,A))).
	$$
\end{proof}

\begin{proof}[Proof of Theorem \ref{Teo formas}]
	We prove the statement for the sum of unstable exponents. One may repeat the argument for $f^{-1}$ to obtain the claim for stable exponents.
	
	Suppose by contradiction that there is a positive volume set  $Z\subset \mathcal{R} \subset {\T}^d$, (where $\mathcal{R}$ is the set of points satisfying Oseledets' theorem as stated in \ref{oseledets} ) such that for all $x\in Z$ we have
	$\sum_{i=1}^{d_u}\lambda^u_i(f,x)>\sum_{i=1}^{d_u}\lambda^u_i(A)$.

	For $q\in {\N}\setminus \{0\}$ we define the set
	$$
	 Z_q=\left\{x\in Z;\sum_{i=1}^{d_u}\lambda_i^u(f,x)>\sum_{i=1}^{d_u}\lambda_i^u(A)+\log\left(1+\frac{1}{q}\right)\right\}.
	$$
	Since $\bigcup_{q=1}^{\infty}Z_q=Z$, there is a $q$ such that $m(Z_q)>0$. For any $x\in Z_q$ we have:
	$$
	\displaystyle\lim_{n\rightarrow\infty}\frac{1}{n}\log|\jac f^n(x)|_{E^u}|>\sum_{i=1}^{d_u}\lambda_i^u(A)+\log\left(1+\frac{1}{q}\right).
	$$
	So, there is  $n_0$ such that for $n\geq n_0$ we have:
	\begin{align*}
		\frac{1}{n}\log|\jac f^n(x)|_{E^u}|&>\sum_{i=1}^{d_u}\lambda_i^u(A)+\log\left(1+\frac{1}{q}\right) \\
		&>\frac{1}{n}\log e^{n\sum_{i=1}^{d_u}\lambda_i^u(A)}+\frac{1}{n}\log\left(1+\frac{1}{q}\right)^{n}.
	\end{align*}
	So we get:
	$$
	|\jac f^n(x)|_{E^u}|>\left(1+\frac{1}{q}\right)^{n}e^{n\sum_{i=1}^{d_u}\lambda_i^u(A)}.
	$$
	By this fact, for each $n>0$ we define the set
	$$
	Z_{q,n}=\left\{x\in Z_q;|\jac f^k(x)|_{E^u}|>\left(1+\frac{1}{q}\right)^{k}e^{k\sum_{i=1}^{d_u}\lambda_i^u(A)},\,\,\forall\,\,k\geq n \right\}.
	$$
	So, there is $N>0$ with $m(Z_{q,N})>0$.

	Now for any $x\in {\T}^d$, let $B_x$ be a foliated box of $W^u_{f}$ around $x$. By compactness,
	we can take a finite cover $\{B_{x_i}\}_{i=1}^{j}$ of ${\T}^d$. As $W^u_f$ is absolutely continuous \cite{brinpesin} ,
	there is $i$ and $x\in B_{x_i}$ such that $\vol_ W^u(f)(B_{x_i}\cap W^u_{f}(x)\cap Z_{q,N})>0$.

	Consider $B_r(x)\subset W^u_f(x)$ satisfying $ {\vol_{W^u(f)}(B_r(x)\cap Z_{q,N})>0}$.
	By Proposition \ref{Prop vol formas},  there is  $K(\epsilon) \in \mathbb{N}$ such that for all $k \geq K(\epsilon)$:
 \begin{align}\label{eq dim maior 4}
		\vol_{W^u(f)}(f^k(B_r(x)))\leq (1+\epsilon)^k\vol_{W^u(A)}(A^k(B_r(x)))\leq(1+\epsilon)^k e^{k\sum\lambda^u_i(A)}\vol_{W^u(A)}(B_r(x)).
	\end{align}

 On the other hand, for  $k\geq N$,
	 \begin{align}\label{eq dim maior 5}
		\vol_{W^u(f)}({f}^k(B_r(x)))&=\int_{B_r(x)}|\jac f^k|_{E^u}|d\vol_{W^u(f)} \nonumber\\
		&\geq\int_{B_r(x)\cap Z_{q,N}}|\jac f^k|_{E^u}|d\vol_{W^u(f)} \nonumber\\
		&>\int_{B_r(x)\cap Z_{q,N}}\left(1+\frac{1}{q}\right)^{k}e^{k\sum_i\lambda_i^u(A)} d\vol_{W^u(f)} \nonumber\\
		&= (1+\frac{1}{q})^{k}e^{k\sum_{i=1}^{d_u}
\lambda_i^u(A)} \vol_{W^u(f)}(B_r(x)\cap Z_{q,N}).
	\end{align}  
	Taking $\epsilon < \frac{1}{q},$  for large enough $k$ the inequalities  (\ref{eq dim maior 4}) and (\ref{eq dim maior 5}) give a contradiction. 
 %Then there is large enough $k$ such that
%	$$
%	\left(1+\frac{1}{q}\right)^{k}e^{k\sum_i\lambda_i^u(A)} \vol(B_r(x)\cap Z_{q,N}) >(1+\varepsilon)^k e^{k\sum\lambda^u_i(A)}\vol(B_r(x)).
%	$$

\end{proof}

\subsection{Proof of Theorem \ref{teor A}}

We use  Theorem \ref{Teo formas} and the following proposition to conclude the proof of Theorem \ref{teor A}.
\begin{proposition}\label{Prop W form}
	Let $W$ be a foliation in ${\R}^d$ of dimension $d_u$. If there is a ($d-d_u$)-plane $P$ in ${\R}^d$ such that $\angle(T_xW,P)>\alpha>0$,
	for all $x\in {\R}^d$,
	then there is a $d_u$-form $\omega$ which is closed and non-degenerate on $W$.
\end{proposition}

\begin{proof}
	Let $B$ be the orthogonal complement of $P$, there is always a closed and non-degenerate form in $B$, in fact,
	just take the volume form $\omega=dx_1\wedge dx_2\wedge\cdots\wedge dx_{d_u}$ which is non-degenerate in $B$
 and $\omega$ is closed. Note that $\omega$
	is degenerate only in $B^{\bot}=P$ and by hypothesis $\angle(T_xW,P)>\alpha>0$. So, 
	$\omega$ is closed and non-degenerate in $T_xW$ for all $x\in{\R}^d$.
\end{proof}

\begin{corollary}
	If $A\colon {\T}^d\rightarrow{\T}^d$ is a linear partially hyperbolic diffeomorphism and $f$ is a $C^2$ conservative diffeomorphism which
is a small $C^1$-perturbation of $A$, then
	$\displaystyle\sum_{i=1}^{d_u}\lambda_i^u(f,x)\leq\displaystyle \sum_{i=1}^{d_u}\lambda_i^u(A)$ and
	$\displaystyle\sum_{i=1}^{d_s}\lambda_i^s(f,x)\geq\displaystyle\sum_{i=1}^{d_s}
	\lambda_i^s(A)$ for Lebesgue almost everywhere $x\in {\T}^d$.
\end{corollary}

\subsection{Proof of Theorem \ref{teor B}}

From \cite[Proposi\c{c}\~{a}o 3.1.2]{daquilema2014hyperbolic} we have:

\begin{proposition}[\cite{daquilema2014hyperbolic}]
	Let $f\colon {\T}^d\rightarrow{\T}^d$ be a partially hyperbolic diffeomorphism homotopic to a linear Anosov diffeomorphism $A$ such that:
	\begin{itemize}
		\item[a)] Each element of the homotopic path is partially hyperbolic diffeomorphism;
		\item[b)] If $f_1$ and $f_2$ are two elements of homotopic path, then $\dime E^{\sigma}(f_1)=\dime E^{\sigma}(f_2)$, $\sigma\in \{s,c,u\},$
	\end{itemize}
	then there exists a  closed  non-degenerate $d_u-$ form on $W^u(f)$.
\end{proposition}

By proposition above $f$ has a $d_u$ closed and non-degenerate on $W^u$, using the inverse $f^{-1}$ we get a $d_s$ closed and
non-degenerate form on $W^s$. By Theorem \ref{Teo formas} we conclude the first part of Theorem \ref{teor B}. In the next subsection, we complete the proof.

\subsubsection{Maximizing Measures} \label{maximizing}

In this section, we prove the last part of the Theorem \ref{teor B}. We write all proofs for the unstable 
bundle.

First, we recall the definition of topological and metric entropy along an expanding foliation. In \cite{HHW} the authors define the notion of topological and metric entropy along unstable foliation and prove the variational principle. 

Let $f\colon  M\rightarrow M$ be a $C^1-$partially hyperbolic diffeomorphism and $\mu$ is an $f-$invariant probability measure. 
For a partition $\alpha$ of $M$ denote $\alpha_0^{n-1}=\displaystyle\vee_{i=0}^{n-1}f^{-i}\alpha$ and by $\alpha(x)$
the element of $\alpha$ containing $x$. Given $\varepsilon>0$, let $\mathcal{P}=\mathcal{P}_{\varepsilon}$ denote the set
of finite measurable partitions of $M$ whose elements have diameters smaller than or equal to $\varepsilon$. For each
$\alpha \in \mathcal{P}$ we define a partition $\eta$ such that $\eta(x)=\alpha(x)\cap W^u_{loc}(x)$ for each $x \in M$, where $W^u_{loc}(x)$
denotes the local unstable manifold at $x$ whose size is greater than the diameter $\varepsilon$ of $\alpha$. 
Let $\mathcal{P}^u=\mathcal{P}^u_{\varepsilon}$ denote the set of partitions $\eta$ obtained in this way. We define the 
\emph{ conditional entropy of  $\alpha$ given $\eta$ with respect to $\mu$} by
$$
H_{\mu}(\alpha|\eta):=-\int_{M}\log\mu^{\eta}_{x}(\alpha(x))d\mu(x)
$$
where $\mu^{\eta}_{x}$ are  disintegration of $\mu$ along $\eta$.

\begin{definition}
The \emph{conditional entropy of $f$ with respect to a measurable partition $\alpha$ given $\eta\in \mathcal{P}^u$} is defined as
$$
h_{\mu}(f,\alpha|\eta)=\displaystyle\limsup_{n\rightarrow \infty}\frac{1}{n}H_{\mu}(\alpha_{0}^{n-1}|\eta).
$$

The conditional entropy of $f$ given $\eta \in \mathcal{P}^u$ is defined as
$$
h_{\mu}(f|\eta)=\displaystyle\sup_{\alpha\in\mathcal{P}}h_{\mu}(f,\alpha|\eta),
$$

and the unstable metric entropy of $f$ is defined as
$$
h_{\mu}^u(f)=\displaystyle\sup_{\eta\in\mathcal{P}^u}h_{\mu}(f|\eta).
$$
\end{definition}

Now, we go to define the unstable topological entropy.

We denote by $\rho^u$ the metric induced by the Riemannian structure on the unstable manifold and let 
$\rho^u_n(x,y)=\max_{0\leq j\leq n-1}\rho^u(f^j(x),f^j(y))$. Let $W^u(x,\delta)$ be the open ball inside 
$W^u(x)$ centered at $x$ of radius $\delta$ with respect to the metric $\rho^u$. Let $N^u(f,\epsilon,n,x,\delta)$
be the maximal number of points in $\overline{W^u(x,\delta)}$ with pairwise $\rho^u_ n-$distances at least $\epsilon$.
We call such set an $(n,\epsilon)$ \emph{$u-$separated set} of $\overline{W^u(x,\delta)}$.

\begin{definition}
The unstable topological entropy of $f$ on $M$ is defined by
$$
h^{u}_{top}(f)=\displaystyle\lim_{\delta\to 0}\displaystyle\sup_{x\in M}h^u_{top}(f,\overline{W^u(x,\delta)}),
$$
where
$$
h^u_{top}(f,\overline{W^u(x,\delta)})=\displaystyle\lim_{\epsilon\to 0}\displaystyle\limsup_{n\to\infty}\frac{1}{n}\log N^u(f,\epsilon,n,x,\delta).
$$
\end{definition}

Let $\mathcal{M}_f(M)$ and $\mathcal{M}_f^e(M)$ denote the set of all $f-$invariant and ergodic probability measures on $M$ respectively.

\begin{theorem}[\cite{HHW}]
Let $f\colon  M\rightarrow M$ be a $C^1-$partially hyperbolic diffeomorphism. Then
$$
h^u_{top}(f)=\sup\{h^u_{\mu}(f): \mu\in\mathcal{M}_f(M)\}.
$$
Moreover,
$$
h^u_{top}(f)=\sup\{h^u_{\nu}(f): \nu\in\mathcal{M}_f^e(M)\}.
$$
\end{theorem}

Moreover, the authors proved that unstable topological entropy coincides with unstable volume growth.

\begin{theorem}[\cite{HHW}]
$h^u_{top}(f)=\chi_u(f)$.
\end{theorem}

As in the hypotheses of Theorem \ref{teor B} the diffeomorphisms $f$ and $A$ are homotopic, using the Theorem \ref{Teo Saghin} we have
$$
h^u_{top}(f)=\chi_u(f)=\log sp(f_{\ast},u)=\log sp(A_{\ast},u)=\chi_u(A)=h^u_{top}(A)
$$ 
Thus
$$
\sum_{i = 1}^{d_u} \lambda^u_i(f, x)=h^u_m(f)\leq h^u_{top}(f)=h^u_{top}(A)=  \sum_{i = 1}^{d_u} \lambda^u_i(A).
$$
Then we conclude that $h^u_m(f)= h^u_{top}(f)$.
\begin{flushright}
$\square$
\end{flushright}

\subsection{Proof of Theorem \ref{TeoDecmDominada}}

%{\color{red} Let $X\subset {\R}^d$ and $R>0$. Denote by $B_R(X)$ the neighbourhood
%$$
%B_R(X)=\{y\in {\R}^d; ||x-y||<R\,\, \mbox{for\,\,some}\,\, x\in X\}.
%$$}

We remember that  $A\colon {\T}^d\rightarrow {\T}^d$ is a linear Anosov diffeomorphism and  $0<\lambda<1<\gamma$ are rates for  hyperbolicity, \mbox{$d_s =\dime (E^s_A)$} and  $d_u =\dime (E^u_A).$ Let $W^s_A, W^u_A$ be the stable and unstable foliation of $A$ and by $\widetilde{W}^{s}_A$, $\widetilde{W}^{u}_A$ we denote their lifts to the universal cover $\mathbb{R}^d$. These foliations are stable and unstable foliations of $\widetilde{A}$ which is a lift of $A$. We use ``$\sim $" for objects in the universal cover.  The norm and distance on $\mathbb{R}^d$ are the lift of adapted norm and corresponding distance where the hyperbolicity conditions of $A$ are satisfied, see Subsection \ref{nonuniformanosov} before the announcement of Theorem \ref{TeoDecmDominada}.

Similar to Proposition 2.5 and Corollary 2.6 of \cite{hammerlindl2013leaf} , we show:

\begin{proposition}\label{PropD3}
  	Let $f\colon {\T}^d\rightarrow {\T}^d$ be a diffeomorphism which admits a dominated splitting  
	\mbox{$TM=E\oplus F$} with $||Df|_E||\leq\hat{\gamma}<\gamma$ and $m(Df|_F) \geq\hat{\lambda}>\lambda$ , homotopic to $A$ such that  \mbox{$\dime (E) =d_s$} and  $\dime (F) =d_u.$ Suppose that the  distributions  $E$ and $F$ are integrable and denote by  $\mathcal{E}$ and $\mathcal{F}$ their respective tangent foliations, then there is $R>0$ such that
	\begin{itemize}
		\item $\mathcal{E}(x)\subset B_{R}(\widetilde{W}^{s}_A(x))$,
		\item $\mathcal{F}(x)\subset B_{R}(\widetilde{W}^{u}_A(x))$,
	\end{itemize}
  where $B_{R}(\widetilde{W}^{s}_A(x)) \subset \mathbb{R}^d$ is the set of points which are at distance $R$ from $W^s_A(x)$. Use similar definition for $B_{R}(\widetilde{W}^{u}_A(x)).$
\end{proposition}

\begin{corollary}\label{Cor. D2}
 Fixing  $x \in \mathbb{R}^n$, if  $ ||x-y||\rightarrow\infty$ and $y\in \mathcal{E}(x)$,   
	then $\frac{x-y}{||x-y||}\rightarrow \widetilde{E}^{s}_A(x)$ uniformly.
	More precisely, for $\varepsilon>0$ there is $M>0$ such that if $x\in{\R}^d$, $y\in \mathcal{E}(x)$ and $||x-y||>M$, then
	$$
	||\pi_A^{u}(x-y)||<\varepsilon||\pi^{s}_A(x-y)||.
	$$
	where $\pi_A^{s}$ is the  orthogonal projection on the subspace $E^s_A$  along $E^u_A$ and $\pi_A^{u}$ is the projection on the  subspace $E^{u}_A$  along $E^s_A.$   
\end{corollary}

Analogous statements hold for $\mathcal{F}$. The following proposition is a topological remark (See  \cite{hammerlindl2013leaf}.) and comes from the fact that $f$ and $A$ are homotopic and $A$ is not singular. 

\begin{proposition}\label{Prop Hamm 1}
	Let $f\colon {\T}^d\rightarrow {\T}^d$ be a homeomorphism with linearization $A$,
	then for each $k \in{\Z}$ and $C>1$ there is an $M>0$ such that  for all $x,y\in {\R}^d$ and
	$$
	||x-y||>M\,\,\, \Rightarrow\,\,\,\frac{1}{C}<\frac{||\tilde{f}^k(x)-\tilde{f}^k(y)||}{||\tilde{A}^k(x)-\tilde{A}^k(y)||}<C.
	$$
	More generally, for each $k\in {\Z}$, $C>1$ and any linear projection $\pi\colon {\R}^d\rightarrow{\R}^d$,
	with $A-$invariant image, there is $M > 0$ such that, for $x,y\in {\R}^d$, with $||x-y||>M$
	$$
	\frac{1}{C}<\frac{||\pi(\tilde{f}^k(x)-\tilde{f}^k(y))||}{||\pi(\tilde{A}^k(x)-\tilde{A}^k(y))||}<C.
	$$
\end{proposition}

\begin{proof}[Proof of Proposition \ref{PropD3}]
	
To prove the first claim (the second one is similar), it is enough to show that $||\pi^u_A(x-y)||$ is uniformly bounded for all $y\in \mathcal{E}(x)$. 
Let $C>1$ close enough to $1$ such that $C \hat{\gamma} < \gamma$ and $1 < \frac{\gamma}{C}.$ Put $k=1$ and $C$ chosen as above in the Proposition \ref{Prop Hamm 1} and get appropriate $M$. By contradiction, suppose  that $||\pi^u_A(x-y)||$ is not bounded. So, there is $y\in \mathcal{E}(x)$ with $||\pi^u_A(x-y)||>M$ and consequently, 
	\begin{align*}
||\pi^u_A(\tilde{f}(x)-\tilde{f}(y))|| > \frac{1}{C}||\pi^u_A(\tilde{A}(x-y))||
   = \frac{1}{C} ||(\tilde{A}( \pi^u_A(x-y))|| \geq \frac{\gamma}{C} \|\pi^u_A(x-y) \| > M
	\end{align*}
which implies that we can use induction and for any $n \geq 1$ obtain:

	$$
	||\pi^u_A(\tilde{f}^{n}(x)-\tilde{f}^{n}(y))||>\frac{\gamma^{n}}{C^n} M.
	$$
	Finally, there is a constant $\eta > 0$ such that:
	\begin{equation}\label{eq021}
	    ||\tilde{f}^{n}(x)-\tilde{f}^{n}(y)||>  \eta \frac{\gamma^{n}}{C^n} M.
	\end{equation}
Now consider a smooth curve $\alpha\colon [a,b]\rightarrow \mathcal{E}(x)$ with $\alpha(a)=x$ and $\alpha(b)=y$ whose length is 
$d_{\mathcal{E}}(x,y)$. Then: 

\begin{align*}
||\tilde{f}^{n}(x)-\tilde{f}^{n}(y)||&\leq d_{\mathcal{E}}(\tilde{f}^{n}(x),\tilde{f}^{n}(y))\leq l(\tilde{f}^{n}(\alpha(t)))
=\int_a^b||\frac{d}{dt}(\tilde{f}^{n}(\alpha(t))||dt\\
&\leq \int_a^b||D\tilde{f}^{n}|_{W^{cs}}||\cdot||(\alpha'(t))||dt< \int_a^b\hat{\gamma}^{n} \|\alpha'(t)||dt,
\end{align*}
implies that
\begin{equation}\label{eq022}
    ||\tilde{f}^{n}(x)-\tilde{f}^{n}(y)||<\hat{\gamma}^{n}d_{\mathcal{E}}(x,y)
\end{equation}

 The equations (\ref{eq021}) and (\ref{eq022}) and $\frac{\gamma}{C} > \hat{\gamma}$ give us a contradiction when  $n$ is large enough.
\end{proof}

\begin{proposition}\label{Prop DD}
	Let $f\colon  M\rightarrow M$ be a $C^2$ conservative non-uniformly Anosov diffeomorphism with $TM = E \oplus F$. If the distribution $E$ is integrable, then the respective foliation $\mathcal{E}$ is upper leafwise absolutely continuous.
\end{proposition}

%$\mathcal{R}=\displaystyle\bigcup_{\varepsilon>0}\displaystyle\bigcup_{l=1}^{\infty}\mathcal{R}_{\varepsilon}^{l}$ the union of Pesin's set from the of 
% about  nonuniformly hyperbolic systems 
%with decomposition $E^{cs}\oplus E^{cu}$ 
%(see, \cite{BarreiraPesin2002lyapunov}). 
%Where $\mathcal{R}_{\varepsilon}^{l}$ is the level set $l$. We remember that for each
%$x\in \mathcal{R}$ there are local manifolds $W^{cs}_r(x)$ 
%and $W^{cu}_r(x)$ tangent to $E^{cs}$ and %$E^{cu}$, respectively. It is known 
%that the length $r$ of manifolds is %bounded by below in each level set 
%$\mathcal{R}_{\varepsilon}^{l}$. Furthermore the family of these manifolds are leafwise absolutely continuous.

Let $\mathcal{R}$ be the regular set in Pesin sense such that $m(\mathcal{R}) = 1$ and $\mathcal{R} = \bigcup_{l=1}^{\infty} \mathcal{R}_l$ is the union of Pesin's block where $m(\mathcal{R}_l) \rightarrow 1, \mathcal{R}_l \subset \mathcal{R}_{l+1}.$ Moreover, the size of Pesin stable manifolds of points in $\mathcal{R}_l$ is larger than some positive constant $r_l.$ In general $r_l \rightarrow 0$ when $l \rightarrow \infty.$ We use an absolute continuity result due to Pesin  \cite{BarreiraPesin2002lyapunov}(Theorem 4.3.1) which essentially shows the absolute continuity of Pesin's stable laminations in each block.

%We call $cs-$block, denoted by $\mathbb{B}_{\varepsilon}^{l}$, a $W_r^{cs}-$foliated box in $\mathcal{R}_{\varepsilon}^{l}$.

\begin{proof}[Proof of Proposition \ref{Prop DD}] 
	Suppose by contraction that $\mathcal{E}$ is not upper leafwise absolutely continuous. So, 
	there exist a foliated box $B$ and a set $Z\subset B$ with $\vol_{\mathcal{E}}(Z\cap \mathcal{E}(x))=0$ 
	for $m-$a.e. $x \in B$ and $m(Z)>0$.
 There exists $l \geq 1$ such that $m(Z \cap \mathcal{R}_l) > 0.$	Take a Lebesgue density point $z\in Z\cap\mathcal{R}_l$ such that for  $V$ a small neighbourhood of $z$ we have $m(Z \cap \mathcal{R}_l \cap V) > 0.$
 
 As Pesin stable manifolds are contained in the leaves of $\mathcal{E}$, our assumption $\vol_{\mathcal{E}}(Z\cap \mathcal{E}(x))=0$ and upper leafwise absolute continuity of stable manifolds  imply that $Z$ has zero measure with respect to conditional measures of $m$ along Pesin stable manifolds  in $\mathcal{R}_l$.)  Hence, $m(Z \cap \mathcal{R}_l \cap V) = 0$ which is a contradiction.
 
 %we conclude that 
 %There is $\varepsilon$, $l$ and a $cs-$block 
%	$\mathbb{B}_{\varepsilon}^{l}$ containing $z$ such that $m(\mathbb{B}_{\varepsilon}^{l}\cap Z)>0$. 
%	Since each Pesin's manifold  $W_r^{cs}(x)$ is containing in the leaf $W^{cs}(x)$, then 
%	$\vol_{W_r^{cs}}(Z\cap W_r^{cs}(x))=0$ for 
%	$m-$a.e. $x\in \mathbb{B}_{\varepsilon}^{l}\cap Z$, but this contradicts the fact of the Pesin's manifolds be upper leafwise absolutely continuous.

\end{proof}

%{\bf Question:} Is $W^{cs}$ lower leafwise absolutely continuous too?

Fix $x\in \mathbb{R}^n$ and denote by $U_r\subset\widetilde{W}^u_A (x)$ the ball of radius $r$ with center $x$.

\begin{proposition}\label{Prop vol Dominada}
	Let $f\colon {\T}^d\rightarrow {\T}^d$ as in the Theorem \ref{TeoDecmDominada}, then given $\varepsilon>0$, there is $r_0>0$ and a constant $C_0>0$ such that
	$$
\vol_{\mathcal{F}}(\tilde{f}^n\pi_x^{-1}(U_r))\leq C_0(1+\varepsilon)^{nd_u}e^{n\sum_{i=1}^{d_u}  \lambda_i^u(A)}\vol_{\tilde{W}_A^u}(U_r)
	$$
	for all $n>0$ and $r \geq r_0$, where $\pi_z$ is the orthogonal projection from $\mathcal{F}(z)$ to $\widetilde{W}_A^u(z)$ 
 (along $\widetilde{E}_A^s$) 
 for any $z \in \mathbb{R}^n$. 
\end{proposition}

%To prove this proposition we need the following auxiliary results:

%\begin{lemma}\label{lema Lp}
%Let  $f:M\rightarrow M$ be a Lipschitzian map, then there is $K > 0$ such that $|\jac f(x)|\leq K,$ for any $x\in M$ such that  $f$ is differentiable at $x.$
%\end{lemma}

%\begin{proof}[Proof of Lemma]
%Let $L$ be the the Lipschitz constant of $f$, give  $x\in M$ such that $f$ is differentiable and  $v\in T_xM$, then
%$$
%||D_xf(v)||=\displaystyle\lim_{t\rightarrow 0}\left\|\frac{f(x + tv)-f(x)}{t}\right\|\leq\lim_{t\rightarrow 0}\frac{L||x+tv-x||}{|t|}=L||v||,
%$$
%it implies that  $||D_xf||=\displaystyle\sup_{v\in T_xM}\frac{||D_xf(v)||}{||v||}\leq L.$

%Since $\R^{n^2}$  is isomorphic to $M_n(\R),$ the space of all $n\times n-$matrices with real coefficients,
%by equivalence between norms of $\R^{n^2}$ we have \mbox{$\{A \in M_n(\R)|\; ||A||\leq L\}$} is compact.
%Since $\det: M_n(\R) \rightarrow \R $ is continuous, then there is $K \geq 0$ such that
%$$
%|\jac f(x)|=|\det D_xf| \leq K.
%$$
%It concludes the proof.
%\end{proof}

%\begin{proposition}\label{Prop recob}
%Let $f:X\rightarrow Y$ be a covering map such that $X$ is path connected  and $Y$ is simply connected, then $f$ is a homeomorphism.
%\end{proposition}

%The proof of the Proposition \ref{Prop vol Dominada} is totally analogous to proof 
%of (Proposition 3.2, \cite{CostaMicena2017pathological}). Let's write it for completeness of the work. 

\begin{proof}[Proof of the Proposition \ref{Prop vol Dominada}]
We prove the following claims.

\textbf{Claim 1:}
 For  $z\in{\R}^d$, the orthogonal  projection \mbox{$\pi_z\colon  \mathcal{F}(z)\rightarrow \widetilde{W}^u_A(z)$}
%(parallel to  $\widetilde{E}^s_A$) 
is a uniform bi-Lipschitz diffeomorphism.
\begin{proof}[Proof of Claim 1]
By item (2) of the Theorem \ref{TeoDecmDominada} and continuity of $E$ and $F$ there is $\alpha>0$ such that the angle $\angle(E,E^u_A)>\alpha$\,\, and
   \,\,$\angle(F, E^s_A)>\alpha.$ So  $\widetilde{W}_A^s$ is uniformly transversal to
the foliation $\mathcal{F}$ and there is  $\beta>0$ such that
		\begin{equation}\label{eq10}
		||d\pi_z(x)(v)||\geq\beta||v||
		\end{equation}
for any $x\in \mathcal{F}(z)$ and $v\in T_x\mathcal{F}(z)$. This implies that $\jac \pi_z(x)\neq 0$.
By Inverse Function Theorem, for each $x\in \mathcal{F}(z)$ there is a ball
$B(x, \delta)\subset \mathcal{F}(z)$ such that $\pi_z|_{B(x, \delta)}$ is a diffeomorphism. In fact from the proof of the Inverse
Function Theorem and Equation (\ref{eq10}), $\delta$ can be taken independent of $x$. Again, by  Equation (\ref{eq10})
there is $\varepsilon>0$, independent of $x$ such that $B(\pi_z(x), \varepsilon)\subset \pi_z(B(x, \delta))$. In order to prove 
 that $\pi_z$ is surjective, we show that $\pi_z(\mathcal{F}(z))$ is an open and closed subset. As $\pi_z$ is a local homeomorphism, then
 $\pi_z(\mathcal{F}(z))$ is open. To verify that it is also closed let
$y_n\in \pi_z(\mathcal{F}(z))$ be a sequence converging to $y$, hence there is  $n_0$ large enough
such that $y\in B(\varepsilon,y_{n_0})$ and therefore $y\in {\pi_z}(\mathcal{F}(z))$. So, $\pi_z$ is surjective.
Moreover, $\pi_z$ is a covering map and injectivity follows from the fact that any covering map from a path-connected space to a simply connected space is a homeomorphism.

Let us prove that $\pi_z$ is bi-Lipschitz.  In fact, $$||\pi_z(x)-\pi_z(y)||\leq||x-y||\leq d_{\mathcal{F}}(x,y).$$ Let us show that
$\pi_z^{-1}$ is also Lipschitz. This is immediate from Equation (\ref{eq10}).
 Indeed, let  $[x,y]$ be the line segment in  $\widetilde{W}^u_A (z)$ connecting $x$ to  $y,$ so the set
$\pi_z^{-1}([x,y])=\Gamma$ is a smooth curve connecting the points $\pi_z^{-1}(x)$ and $\pi_z^{-1}(y)$ in $\mathcal{F}(z)$. So,
$$
d_{\mathcal{F}}(\pi_z^{-1}(x),\pi_z^{-1}(y))\leq length(\Gamma)=\int_{[x,y]}|d\pi_z^{-1}(t)|dt\leq \frac{1}{\beta}||x-y||.
$$
\end{proof}

%As $E^u(A)=E_1\oplus E_2\oplus\cdots\oplus E_{d_u}$, suppose that the edges of the  $d_u$-cube $R$ lies in hypercubes parallel to  $E_i.$ \\

%As in  (Lemma 3.6, \cite{MT}), we can prove:

%Denote by  $\nu=\min\{|\lambda_i^u(A)|\}$ and $\eta=\max\{|\lambda_i^u(A)|\}$, where $\lambda_i^u(A)$ are the eigenvalues of $A|_{E^u_A}$ (possibly complex). If the eigenvalues is complex, $|\lambda_i^u(A)|$ represents its norm. 

%\textbf{Claim 2:}  Given  $\varepsilon>0$, there is $M>0$ such that if $y\in \widetilde{E}^u_A(x)$ and $||x-y||>M$ then;
%$$(1+\varepsilon)^{-1}||\tilde{A}(x)-\tilde{A}(y)||\leq  ||\tilde{A}\pi_x^{-1}(x)-\tilde{A}\pi_x^{-1}(y)|| \leq (1+\varepsilon)||\tilde{A}(x) -\tilde{A}(y)||.
%$$
%\begin{proof}[Proof of Claim 2]
%Observe that \begin{equation} \label{ineq}
%||\tilde{A}\pi_x^{-1}(x)-\tilde{A}\pi_x^{-1}(y)|| \leq ||\tilde{A}(x)-\tilde{A}(y)|| + ||\tilde{A}y-\tilde{A}\pi_x^{-1}(y)||.\end{equation}
%Now, as $M$ is large enough, by Corollary \ref{Cor. D2}, $\|y -\pi_x^{-1}(y)\|$ is much smaller than $\|x-y\|$ and $||\tilde{A}y-\tilde{A}\pi_x^{-1}(y)|| \leq \varepsilon ||\tilde{A}(x)-\tilde{A}(y)||.$ Substituting in \ref{ineq} we obtain the proof of (one part) the claim. Another inequality in the claim is proved similarly.
%\end{proof}

\begin{figure}[!htb]
\centering
\includegraphics[scale=0.75]{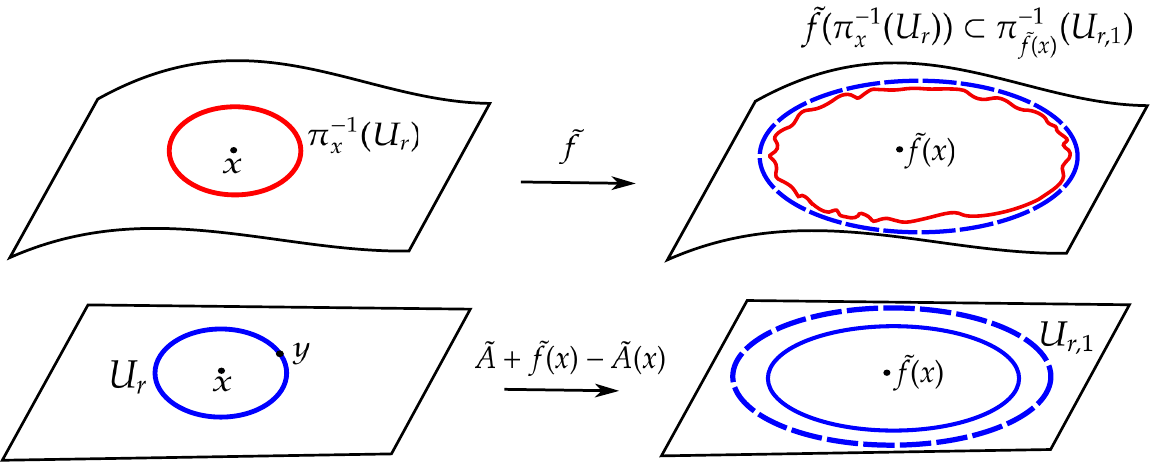}
\caption{}\label{cone1}
\end{figure}

%\textbf{Claim 3:}  There is $M > 0,$ such that $\tilde{f}^n(\pi^{-1}(R)) \subset \pi^{-1}((1+ \varepsilon)^n\tilde{A}^n(R)),$ for every $n \geq 1,$ up to translation of $(1+ \varepsilon)^n\tilde{A}^n(R).$  We are considering $R$  a $d_u$ hypercube with edges parallel to each direction $E^u_i, i =1, \ldots, d_u,$ (see Figure 1), and the length of each edge of $R$ is at least $M.$

\textbf{Claim 2:} Given $\varepsilon>0$, there is $r_0 > 0$ such that $\tilde{f}^n(\pi_x^{-1}(U_r)) \subset \pi_{\tilde{f}^n(x)}^{-1}(U_{r, n}),$ for every $r \geq r_0$ and  $n \geq 1,$ where $U_{r, n} \subset  \widetilde{W}_A^u (\tilde{f}^n(x))$ is a $d_u-$dimensional ellipsoid  with volume bounded above by   $(1+\epsilon)^{nd_u} \vol_{\tilde{W}_A^u}(A^n(U_r)).$
\begin{proof}
%Take $\varepsilon>0$  small enough such that $(1+\varepsilon)^{-2} m(\tilde{A}|_{E^u}) > 1$ and  $C:=1+\varepsilon$. By the Proposition \ref{Prop Hamm 1}, there is $M>0$ such that if $y\in \widetilde{E}^{u}_A (x)$ with 
%$ ||x - y|| \geq M$ we have
%\begin{equation}\label{eq023}
 %   ||\tilde{f}(x)-\tilde{f}\pi_x^{-1}(y)||\leq (1+\varepsilon)||
%\tilde{A}(x)-\tilde{A}\pi_x^{-1}(y)||.
%\end{equation}

Let $ r \geq r_0 :=  \frac{2R+2K}{\varepsilon m(A|_{E^u})}$ where $R$ comes from Proposition \ref{PropD3} and \mbox{$\|\tilde{f}- \tilde{A}\|_{\infty} \leq K.$}   We prove that \begin{equation} \label{inclusion}
\tilde{f}(\pi_x^{-1}(U_r)) \subset \pi_{\tilde{f}(x)}^{-1}(U_{r, 1}), 
\end{equation}
where $U_{r,1} \subset \widetilde{W}^u_A(\tilde{f}(x))$ is obtained from $\tilde{A}(U_r)$ by first applying a homothety of ratio $(1+\epsilon)$ centered at $\tilde{A}(x)$ and then translating by $\tilde{f}(x) - \tilde{A}(x)$.  Observe that $U_{r, 1}$ is an ellipsoid inside the affine $d_u$ dimensional subspace passing through $\tilde{f}(x).$
 Take any $y$ on the boundary of $U_r$. let $z:=\pi_{\tilde{f}(x)}(\tilde{f}(\pi_x^{-1}(y))) \in \widetilde{W}^u_A(\tilde{f}(x))$.
On the one hand, we have 
\begin{align*}  \| z - (\tilde{A}(y) + \tilde{f}(x)- \tilde{A}(x))   \| & \leq \| \tilde{f}(x)- \tilde{A}(x)   \| + \| z - \tilde{A}(y)  \| 
& \numberthis \label{first} \leq K + \| z - \tilde{A}(y) \|.  \end{align*}
On the other hand,
\begin{align*}
  \| z - \tilde{A}(y)  \| \leq \|z - \tilde{f}(\pi_x^{-1}(y))\| &+ \|  \tilde{f}(\pi_x^{-1}(y))-\tilde{A}(\pi_x^{-1}(y)) \|  \\
  &+ \|   \tilde{A}(\pi_x^{-1}(y)) - \tilde{A}(y) \| \\
  & \numberthis \label{second} \leq 2R+K. \end{align*}
  
In the above inequalities we have used two times Proposition \ref{PropD3} to get $$\|z - \tilde{f}(\pi_x^{-1}(y)\| \leq R, $$ $$\|   \tilde{A}(\pi_x^{-1}(y)) - \tilde{A}(y) \| \leq \|y - \pi_x^{-1}(y)\| \leq R, $$ (observe that $y - \pi_x^{-1}(y)$ belongs to the stable subspace of $A$) and finally $$ \| \tilde{f}(\pi_x^{-1}(y)))-\tilde{A}(\pi_x^{-1}(y)) \| \leq \|\tilde{f}-\tilde{A}\|_{\infty} \leq K.$$ 

Now, putting  inequalities (\ref{first}) and (\ref{second}) together we get:
$$ \|z - (\tilde{A}(y) + \tilde{f}(x)- \tilde{A}(x))\| \leq 2K + 2R.$$ Observe that $\tilde{A}(y) + \tilde{f}(x)- \tilde{A}(x)$ is the translation of the $\tilde{A}(y)$ and belongs to $\widetilde{W}_A^u(\tilde{f}(x)).$  Indeed,  $\tilde{A}(y)- \tilde{A}(x)$ is a vector which belongs to the unstable bundle of $\tilde{A}$ and we identify the unstable bundle at $\tilde{f}(x)$ with the corresponding affine subspace of $\mathbb{R}^n$ passing through $\tilde{f}(x)$. 
Now, as the distance between $(1+\epsilon) \tilde{A}(U_r) $ and $ \tilde{A}(U_r)$ is $ \epsilon r m(\tilde{A}|_{E^u}) \geq 2R + 2K$ (by the choice of $r$), 
we conclude that $z$ belongs to the   ellipsoid $U_{r, 1} :=(\tilde{f}(x)- \tilde{A}(x)) + (1+\epsilon) \tilde{A}(U_r)$ which proves  (\ref{inclusion}). Moreover, observe that $\vol_{\widetilde{W}^u_A} (U_{r, 1}) =(1+\epsilon)^{d_u} \vol_{\widetilde{W}^u_A}(\tilde{A}(U_r)) = (1+\epsilon)^{d_u} e^{\sum_{i=1}^{d_u} \lambda^u_i(A)} \vol_{\widetilde{W}^u_A} (U_r).$
%{\color{} Essa ultima desigualdade precisa de constant, ou talvez assumirmos queue todos os fibrados dentro de Eu sao ortogonais!?}

%It is important to note that if $||x - y|| \geq M$, then by the choice of $\varepsilon$,
%$$||\tilde{f}(x)-\tilde{f}\pi_x^{-1}(y)||>M.$$ Indeed, by Proposition \ref{Prop Hamm 1} and  Claim 2, we have
%\begin{align*}
%||\tilde{f}(x)-\tilde{f}\pi_x^{-1}(y)||&\geq (1+\varepsilon)^{-2}||
%\tilde{A}(x)-\tilde{A}(y)||\\
%& (1+\varepsilon)^{-2} m(A|_{E^u}) M \geq M.
%\end{align*}
Now we apply again $\tilde{f}$ and obtain $$\tilde{f}^2(\pi_x^{-1}(U_r)) \subset \tilde{f} (\pi_{\tilde{f}(x)}^{-1}(U_{r,1})) $$
As the distance between $\tilde{A}(U_{r, 1})$ and $(1+\epsilon)\tilde{A}(U_{r, 1})$ is larger than $\varepsilon r m(\tilde{A}|_{E^u})$, similarly as above we obtain
$$ \tilde{f} (\pi_{\tilde{f}(x)}^{-1}(U_{r,1})) \subset \pi_{\tilde{f}^2(x)}^{-1}(U_{r, 2})$$ where $U_{r, 2}$ is a translation of $(1+\varepsilon)^2 \tilde{A}^2(U_r).$
 In fact, inductively we obtain:
 \begin{equation} \label{inclusions} \tilde{f}^{n+1}(\pi_x^{-1}(U_r))  \subset \tilde{f} (\tilde{f}^{n}(\pi_x^{-1}(U_r))) \subset \tilde{f} (\pi_{\tilde{f}^n(x)}^{-1}(U_{r, n})) \subset \pi_{\tilde{f}^{n+1}(x)}^{-1}(U_{r, n+1}) 
 \end{equation}
where $U_{r, n+1}$ is an ellipsoid with volume less than $(1+\varepsilon)^{(n+1)d_u} e^{(n+1) \sum_{i=1}^{d_u} \lambda^u_i(A)} \vol_{\tilde{W}_A^u}(U_r).$ The last inclusion in (\ref{inclusions}) follows using the same arguments as above to prove  (\ref{inclusion}) substituting the ball $U_r$ by the ellipsoid $U_{r, n}.$
\end{proof}

By Claim 1, for all $x\in M$ we have $|\jac\, \pi^{-1}_x|$ is uniformly bounded and consequently  there is a constant $C_0 >0$ such that
\begin{align*}
\vol_{\mathcal{F}}(\tilde{f}^n\pi_x^{-1}(U_r))&\leq\vol_{\mathcal{F}}(\pi_{\tilde{f}^n(x)}^{-1}(U_{r, n})  \leq C_0(1+ \varepsilon)^{nd_u} e^{n \sum_{i = 1}^{d_u} \lambda^u_i(A)}\vol_{\tilde{W}_A^u}(U_r).
\end{align*}
 It concludes the proof of the Proposition \ref{Prop vol Dominada}.
\end{proof}

\begin{proof}[Proof of Theorem \ref{TeoDecmDominada}]

Suppose by contradiction that there is a positive volume set  $Z\subset {\T}^d$ such that
$\sum_{i=1}^{d_u}\lambda_i^{cu}(f,x)>\sum_{i=1}^{d_u}\lambda_i^u(A)$ for any $x\in Z$.

Let $P\colon {\R}^d\rightarrow {\T}^d$ be the covering map, $D\subset{\R}^d$ a fundamental domain  and $\widetilde{Z}=P^{-1}(Z)\cap D$.  We have $\vol(\widetilde{Z})>0$.

For each $q\in {\N}\setminus \{0\}$ we define the set
$$
Z_q=\left\{x\in \widetilde{Z};\sum_{i=1}^{d_u}\lambda_i^{cu}(\tilde{f},x)>\sum_{i=1}^{d_u}\lambda_i^u(A)+\log\left(1+\frac{1}{q}\right)\right\}.
$$
We have $\bigcup_{q=1}^{\infty}Z_q=\widetilde{Z}$, thus there is  $q$ such that $m(Z_q)>0$. For each $x\in Z_q$ it follows that
$$
\displaystyle\lim_{n\rightarrow\infty}\frac{1}{n}\log|\jac \tilde{f}^n(x)|_{ F}|>\sum_{i=1}^{d_u}\lambda_i^u(A)+\log\left(1+\frac{1}{q}\right).
$$
So there is  $n_0$ such that for   $n\geq n_0$  we have
\begin{align*}
\frac{1}{n}\log|\jac \tilde{f}^n(x)|_{F}|&>\sum_{i=1}^{d_u}\lambda_i^u(A)+\log\left(1+\frac{1}{q}\right)\\
                                         &>\frac{1}{n}\log e^{n\sum_{i=1}^{d_u}\lambda_i^u(A)}+\frac{1}{n}\log\left(1+\frac{1}{q}\right)^{n}.
\end{align*}
This implies that
$$
|\jac \tilde{f}^n(x)|_F|>\left(1+\frac{1}{q}\right)^{n}e^{n\sum_{i=1}^{d_u}\lambda_i^u(A)}.
$$
For every $n>0$ we define
$$
Z_{q,n}=\left\{x\in Z_q;|\jac \tilde{f}^k(x)|_F|>\left(1+\frac{1}{q}\right)^{k}e^{k\sum_{i=1}^{d_u}\lambda_i^{u}(A)},\,\,\forall\,\,k\geq n \right\}.
$$
There is $N>0$ with $\vol(Z_{q,N})>0$.

 For each $x\in D$, consider $B_x\subset{\R}^d$ a foliated box of  $\mathcal{F}.$
By compactness there is finite cover $\{B_{x_i}\}_{i=1}^{j}$ covering $\overline{D}.$ Since $W^{cu}_f$ is absolutely continuous and the covering map is smooth, then  $\mathcal{F}$ is absolutely continuous, thus there is some $i$ and $p\in B_{x_i}$
such that $\vol_{\mathcal{F}}(B_{x_i}\cap \mathcal{F}(p)\cap Z_{q,N})>0$. 

There is a set $\pi_p^{-1}(U_r)\subset \mathcal{F}(p)$, where $\pi_p^{-1}(U_r)$ is as in the Proposition \ref{Prop vol Dominada}
containing $p$ and $r$ is large enough 
such that \mbox{$\vol_{\mathcal{F}}(\pi_p^{-1}(U_r)\cap Z_{q,N})>0.$}
Let  $\alpha>0$ be such that $\vol_{\mathcal{F}}(\pi_p^{-1}(U_r)\cap Z_{q,N})=\alpha \vol_{\mathcal{F}}(\pi_p^{-1}(U_r)),$
by Proposition \ref{Prop vol Dominada} we have

\begin{equation}\label{eq4}
\vol_{\mathcal{F}}(\tilde{f}^n(\pi_p^{-1}(U_r)))\leq C(1+\varepsilon)^de^{n\sum\lambda_i^u(A)}\vol_{\tilde{E}_A^u}(U_r).
\end{equation}
On the other hand,
\begin{align}\label{eq5}
\vol_{\mathcal{F}}(\tilde{f}^n(\pi_p^{-1}(U_r)))&=\displaystyle\int_{\pi_p^{-1}(U_r)}|\jac \tilde{f}^n(x)|_{F}| d \vol_{\mathcal{F}} \nonumber\\
                &\geq\displaystyle\int_{\pi_p^{-1}(U_r)\cap Z_{q,N}}|\jac \tilde{f}^n(x)|_{F}| d \vol_{\mathcal{F}}\nonumber\\
                &>\displaystyle\int_{\pi_p^{-1}(U_r)\cap Z_{q,N}}\left(1+\frac{1}{q}\right)^{n}e^{n\sum_i\lambda_i^u(A)}d \vol_{\mathcal{F}}\nonumber\\
                &>\left(1+\frac{1}{q}\right)^{n}e^{n\sum_i\lambda_i^u(A)} \vol_{\mathcal{F}}(\pi_p^{-1}(U_r)\cap Z_{q,N})\nonumber\\
                &>\left(1+\frac{1}{q}\right)^{n}e^{n\sum_i\lambda_i^u(A)} \alpha\vol_{\mathcal{F}}(\pi_p^{-1}(U_r))\nonumber\\
                &>\left(1+\frac{1}{q}\right)^{n}e^{n\sum_i\lambda_i^u(A)} \alpha\vol_{\tilde{E}_A^u}(U_r)
\end{align}
The equations (\ref{eq4}) and (\ref{eq5}) give us a contradiction when  $n$ is large enough, thus proving the Theorem \ref{TeoDecmDominada}.

\end{proof}

%{\color{red}

%\textbf{Acknowledgements}:  J.S. Costa would like to thank A. Tahzibi for doctoral supervision and several conversations and also to thank Mario R. R. Daquilema for valuable conversations. J. S. Costa was supported by CAPES-PROEX and CNPq process 141224/2013-4. A. Tahzibi was supported by FAPESP thematic project 2017/06463-3 and CNPq productivity fellowship.

%}

\newpage

%--------------------------------------------------------------------------------------

\bibliographystyle{abbrv}

\bibliography{referencesCostaTahzibi}

\end{document}